\documentclass[a4paper,12pt,reqno]{amsart}
\usepackage[T1]{fontenc}
\usepackage[utf8]{inputenc}
\usepackage{amssymb,dsfont,mathrsfs,mathabx}
\usepackage[margin=1in]{geometry}
\usepackage{enumitem}
\usepackage[bookmarksdepth=2]{hyperref}
\usepackage[alphabetic, initials]{amsrefs}
\usepackage{graphicx,xcolor}

\numberwithin{equation}{section}

\theoremstyle{definition}

\theoremstyle{notation}

\newcommand{\R}{\mathds{R}}

\allowdisplaybreaks[1]

\title[A simple but effective bushfire model]{A simple but effective bushfire model: \\
analysis and real-time simulations}

\author[S. Dipierro]{Serena Dipierro}
\address{S. D., 
Department of Mathematics and Statistics,
University of Western Australia,
35~Stirling Highway, WA 6009 Crawley, Australia. }

\email{serena.dipierro@uwa.edu.au}

\author[E. Valdinoci]{Enrico Valdinoci}
\address{E. V., 
Department of Mathematics and Statistics,
University of Western Australia,
35~Stirling Highway, WA 6009 Crawley, Australia. }

\email{enrico.valdinoci@uwa.edu.au}

\author[G. Wheeler]{Glen Wheeler}
\address{G. W.,
School of Mathematics and Applied Statistics,
University of Wollongong,
Northfields Avenue, NSW 2500 Wollongong, Australia. }

\email{glenw@uow.edu.au}

\author[V.-M. Wheeler]{Valentina-Mira Wheeler}
\address{V.-M. W.,
School of Mathematics and Applied Statistics,
University of Wollongong,
Northfields Avenue, NSW 2500 Wollongong, Australia. }

\email{vwheeler@uow.edu.au}

\thanks{Supported by Australian Research Council FL190100081 and DE190100379.}

\begin{document}

\begin{abstract}
We introduce a simple mathematical model for bushfires accounting for temperature diffusion in the presence of a combustion term which is activated above a given ignition state. The model also takes into consideration the effect of the environmental wind and of the pyrogenic flow.

The simplicity of the model is highlighted from the fact that it is described by a single scalar equation, containing only four terms, making it very handy for rapid and effective numerical simulations which run in real-time.

In spite of its simplicity, the model is in agreement with data collected from bushfire experiments in the lab, as well as with spreading of bushfires that have been observed in the real world.

Moreover, the equation describing the temperature evolution can be easily linked to a geometric evolution problem describing the level sets of the ignition state.
\end{abstract}

\maketitle

\section{The bushfire model}

The study of bushfire propagation requires, in general, the understanding of very complex and intertwined phenomena, involving fluid mechanics, physics, chemistry, meteorology, and thermodynamics. In this spirit, a thorough bushfire analysis should combine
weather and climate models with fuel and terrain models, describing accurately flaming
combustion phenomena and eddy formation in different spatial and temporal scales
(see e.g.~\cite{BAK}).
The ``ultimate'' universal bushfire model is not available yet,
and it is customary to approach the problem under different perspectives,
so that often the different models are classified into the categories of
``physical'' (or ``theoretical''), ``semi-physical'' (or ``semi-empirical''), and ``empirical'', see e.g.~\cites{PERRY, PASTOR}; sometimes ``mathematical'' models have been considered as a separate category closely related to ``simulation'', see e.g.~\cites{SULL1, SULL2, SULL3}.

Of course, on the one hand, it is desirable to describe bushfires by a model which is as complete and accurate as possible, but, on the other hand,
it would also be important to have ``simple'' (though not simplistic) models allowing for
rigorous mathematical treatments, a qualitative understanding of the main processes involved,
as well as fast simulations making the model reliable but also concretely usable in real-time.
These simple models
can also act as initial blueprints from which more complex models can be built.

With these considerations in mind, the model that we propose consists of a single scalar equation, in which the effects coming from fluid dynamics and chemical combustion
are encoded by state parameters and interaction kernels. More specifically, the equation on which we focus our attention has the form
\begin{equation}\label{MAIN:EQ}
\partial_t u=c\Delta u+ \int_{\R^n}\big(u(y,t)-\Theta(y,t)\big)_+\,K(x,y)\,dy+
\left( \left(\omega-\frac{\beta(u)\nabla u}{|\nabla u|^\alpha}\right)\cdot\nabla u\right)_-\,,
\end{equation}
where $ u$, $\Theta: \R^n \times \R \rightarrow \R $ and~$ c: \R^n \times \R \rightarrow (0,+\infty) $ are real functions depending on the space point $ x \in \R^n $ and time $ t \in \R $, $ K: \R^n \times \R^n \rightarrow [0,+\infty) $ 
is an interaction kernel,
$ \omega: \R^n \times \R \rightarrow \R^n $ is a time-dependent vector field, $ \beta: \R \rightarrow \R $, $ \alpha \in[0,2] $, and we used the notation for the ``positive and negative parts'' of a function
\begin{equation}\label{ponep}
{\mbox{$f_+(x):=\max\{0,f(x)\}$ and~$f_-(x):=\max\{0,-f(x)\}$.}}\end{equation}

In this model, $u=u(x,t)$ represents the temperature of the environment at the space point~$x\in\R^n$, at time~$t$ (in many applications, one can focus on the case~$n=2$, and deal with a given region with suitable, e.g. Dirichlet, conditions).

As customary, one assumes that temperature spreads through a {\bf heat equation}, in which~$c=c(x,t)$ denotes the environmental diffusion coefficient (the case of constant~$c$ modeling a homogeneous environment).\medskip

The specific model for bushfires accounts for three additional phenomena, namely combustion, wind, and convection.
\medskip

The {\bf combustion} is described by an effective ignition temperature~$\Theta=\Theta(x,t)$, above which the fuel burns.
The case of a constant~$\Theta(x,t)=\Theta_0$ represents the ideal situation of fuel homogeneously spread across the environment, with infinite abundance of fuel at every point.
In this case, the fuel burns at temperature~$\Theta_0$
and affects nearby points through the interaction potential~$K$
(notice that when the temperature~$u$ of the fuel is below the effective ignition temperature~$\Theta$, the combustion effect disappears, in view of the positive part of the function considered in the integral).

Situations in which the availability of fuel varies in space and time are of course modeled by nonconstant functions~$\Theta$. We stress that the model can also include a memory effect, in which the burned fuel turns into ash and is not any longer subject to combustion. This can be done, for instance, by considering the case in which
\begin{equation}\label{TE:d1} \Theta(x,t):=\begin{cases}\displaystyle
\tan\left( \frac{\pi}{2 F} \int_0^ t\big(u(x,\tau)-\bar\Theta(x)\big)_+\,d\tau \right) & {\mbox{if }}\displaystyle\int_0^ t\big(u(x,\tau)-\bar\Theta(x)\big)_+\,d\tau<F,\\
+\infty&{\mbox{otherwise.}}\end{cases}\end{equation}
In this framework, $F=F(x,t)$ accounts for the fuel available
(for instance, $F$ can be taken to be constant;
the case~$F=F(x)$ corresponds to the situation of an initial amount of fuel~$F(x)$ available at the point~$x$, and the case
of~$F$ also depending on time allows one to consider situations in which the fuel is recharged over time).

Moreover, the constitutional ignition temperature at a given point~$x$
is~$\bar\Theta(x)$ (the case of a constant ignition temperature~$\bar\Theta(x)=\Theta_0$ accounting for a uniformly distributed fuel).
Above this ignition temperature, the fuel burns, and the amount of fuel burned is described by the quantity
\begin{equation}\label{TE:d2} B(x,t):= \int_0^ t\big(u(x,\tau)-\bar\Theta(x)\big)_+\,d\tau.\end{equation}

We highlight the difference between the constitutional ignition temperature~$\bar\Theta$ and the effective ignition temperature~$\Theta$. Namely, while~$\bar\Theta$ can be thought as a specific property of the combustible (which is modeled to ignite above this temperature, in the ideal situation of infinite amount of fuel available), the quantity~$\Theta$ takes into account the past state of the fuel and the propellant~$B$ already burned.
Indeed, comparing~\eqref{TE:d1} and~\eqref{TE:d2}, the effective ignition temperature in this case can be written in the form
\begin{equation*} 
\Theta=\tan \left( \frac{\pi B}{2F}\right),\end{equation*}
hence when the quantity of burned combustible~$B$ reaches the threshold~$F$, the fuel gets extinguished, and the corresponding effective ignition temperature diverges to~$+\infty$, making the interaction ignition term in~\eqref{MAIN:EQ} vanish, due to the presence of the positive part in the equation.

See also~\cite{MR2542721}
for models taking into account the environmental temperature with a fuel ignition term.
\medskip

The effect of the environmental {\bf wind} is modeled by
a vector~$\omega=\omega(x,t)$ inducing a transport term in equation~\eqref{MAIN:EQ}.
\medskip

The {\bf convection} is encoded by an additional wind term of the form~$\frac{\beta(u)\nabla u}{|\nabla u|^\alpha}$.
The function~$\beta$ can be thought as supported above the ignition temperature, so that this term is active only in the burning region, directed along, or opposite to, depending on the sign of~$\beta$,
the gradient temperature. The exponent~$\alpha$ can be taken in the interval~$[0,2]$ to further modulate the intensity of this term with respect to the temperature gradient.

The negative part in the term regarding the environmental wind and the convective term accounts for the fact that wind can spread the fire faster, but cannot stop it, hence these terms provide a contribution only when their direction is oriented\footnote{We stress that the total wind term, consisting in the superposition of environmental and pyrogenic winds, namely the term~$\left(\omega-\frac{\beta\nabla u}{|\nabla u|^{\alpha}}\right)\cdot\nabla u$,
provides a contribution only when this function takes negative values. This means that the wind and convective terms provide a contribution only when their direction is oriented opposite to the gradient temperature, that is in the direction of the fire front (which goes from higher to lower values of the temperature, and note that,
in our notation~\eqref{ponep}, the negative part of a function is nonnegative).

In the absence of environmental wind, the pyrogenic flow takes the form~$\frac{(-\beta)_-}{|\nabla u|^{\alpha-2}}=\frac{\beta_+}{|\nabla u|^{\alpha-2}}$, contributing only for positive values of~$\beta$.} opposite to the gradient temperature, which is in turn a proxy for the direction of the front.
In this spirit, the role of the negative part in the last term of~\eqref{MAIN:EQ} is to {\em spread} the fire in the direction of the wind: as the simulations confirm, without the negative part, the burning region will tend to mainly {\em translate} in the direction of the wind.
\medskip

Terms related to a fire-induced wind (also known as ``fire wind'' or ``pyrogenic flow'')
have been introduced in the literature to account for the significant buoyant upflow
over the fire region, created by the strong temperature gradient,
producing local low pressure, which acts as a sink with horizontal pressure gradient
and draws in the surrounding air (see~\cite{SMITH}).

The possibility of modeling the wind advection by a transport term, in various forms,
has also been previously considered in the literature, see e.g.~\cite[equation~(2.2)]{MR4439510}.
In terms of concrete applications, the combination of pyrogenic and environmental wind can be a delicate issue: indeed, on the one hand, in case of wind-driven fires
it is expected the effect of pyrogenic wind to be less prominent (see~\cite{BEER}),
on the other hand pyrogenic winds can be produced even by relatively small wildfires
(as observed in~\cite{LAREAU}).
See~\cites{HILTON1, HILTON2, HILTON201812, MR4439510} and the references therein for different models
describing bushfires accounting for wind and convection effects. 
\medskip

The model that we propose is also prone to modifications and generalizations to capture additional details. For instance, the interaction kernel~$K$ is supposed here to be a nonnegative function (the case in which~$K(x,y)=\bar K(|x-y|)$ describing an ignition interaction
depending only on the distance between sites), but one can replace the term~$K(x,y)\,dy$ by a more general measure integral: for instance, integrating against a Delta Function would return a source term of the type~$\big(u(x,t)-\Theta(x,t)\big)_+$.

In the same spirit, the terms linear in~$(u-\Theta)_+$ can be replaced by nonlinear modifications, to capture nonlinear scalings in fire propagation and radiation effects, the interaction kernels may also depend on the temperature, the Laplacian can be substituted with other diffusive operators, e.g. to model various types of anisotropy and inhomogeneity in the terrain, and
the negative part in the wind term can be dropped, e.g. to model possible containment effects of the wind.\medskip

The Laplacian can also be replaced, or complemented, by nonlinear or nonlocal operators, with the aim of modeling anomalous types of diffusion, and also by operators of mixed orders, for instance
to describe the combination of two types of diffusion (such as the spread of heat along the ground surface and the boundary effect on the ground coming from the spread of heat in the air). For the sake of simplicity, we do not present all these generalizations in full detail here, but we plan to address them specifically in future projects.\medskip

As a sanity check, let us remark that, since both the positive and negative parts of a function are, in our notation~\eqref{ponep}, nonnegative, the solutions of~\eqref{MAIN:EQ} are always {\em supersolutions of the heat equation} (which is consistent, from the physical point of view, with the fact that bushfires can only increase the environmental temperature).\medskip

Albeit non-standard, the partial differential equation
in~\eqref{MAIN:EQ} is of parabolic type and is endowed with a solid mathematical theory, developed in~\cite{PAPER2}.

\section{Numerical validation and real-time simulations}

The model presented in~\eqref{MAIN:EQ} is simple enough to allow for an agile
and effective numerical implementation, producing sufficiently precise results
in real-time.
As concrete examples, we confronted our numerical outputs with supervised
experiments of fires performed in the lab as well as with real data collected
on the occasion of natural bushfires.

\begin{figure}[htbp]
    \centering
    \begin{minipage}[c]{0.49\textwidth}
        \includegraphics[width=\textwidth,trim=20pt 15pt 30pt 25pt,clip]{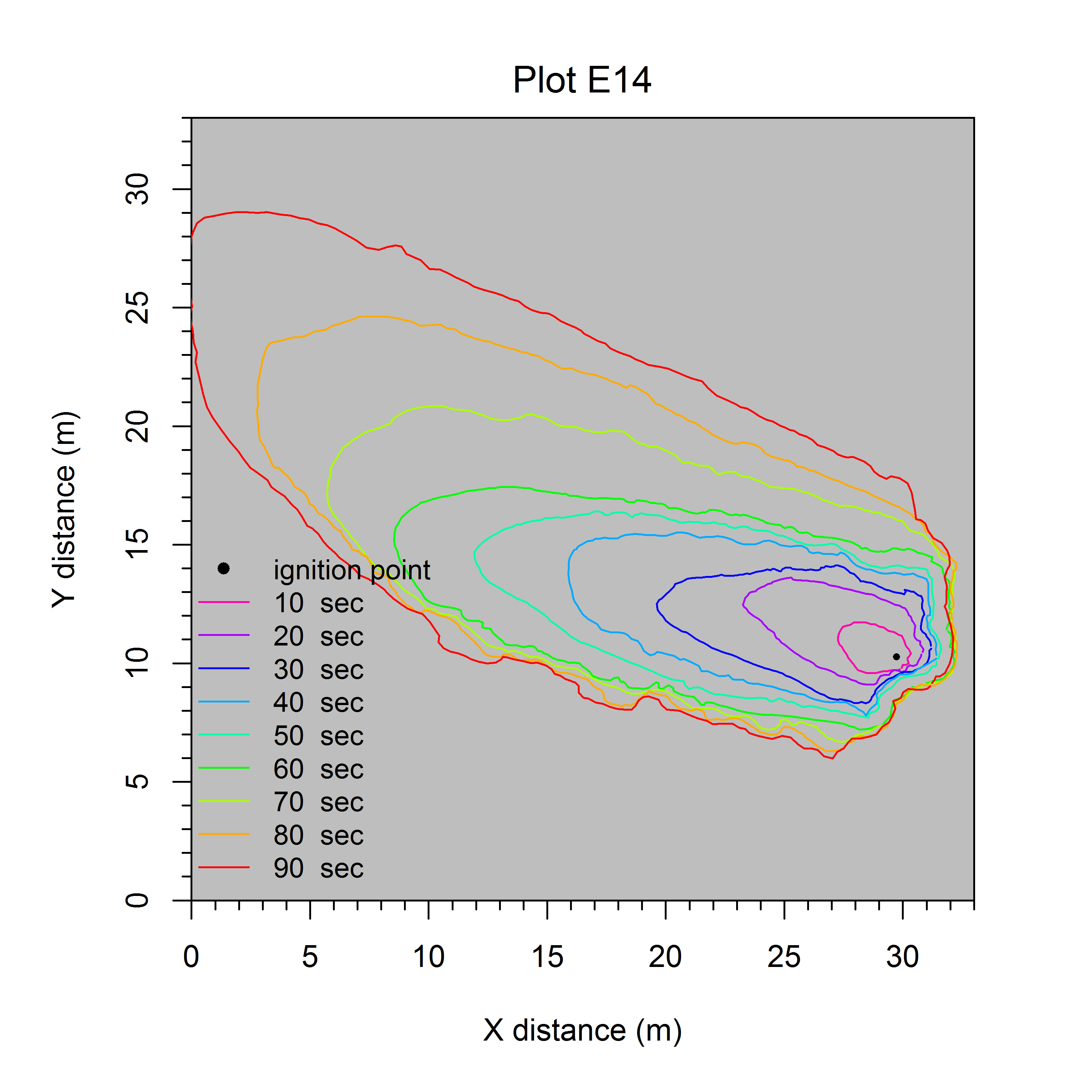}
    \end{minipage}
    \hfill
    \begin{minipage}[c]{0.49\textwidth}
        \centering
        \includegraphics[width=\textwidth,trim=20pt 20pt 30pt 25pt,clip]{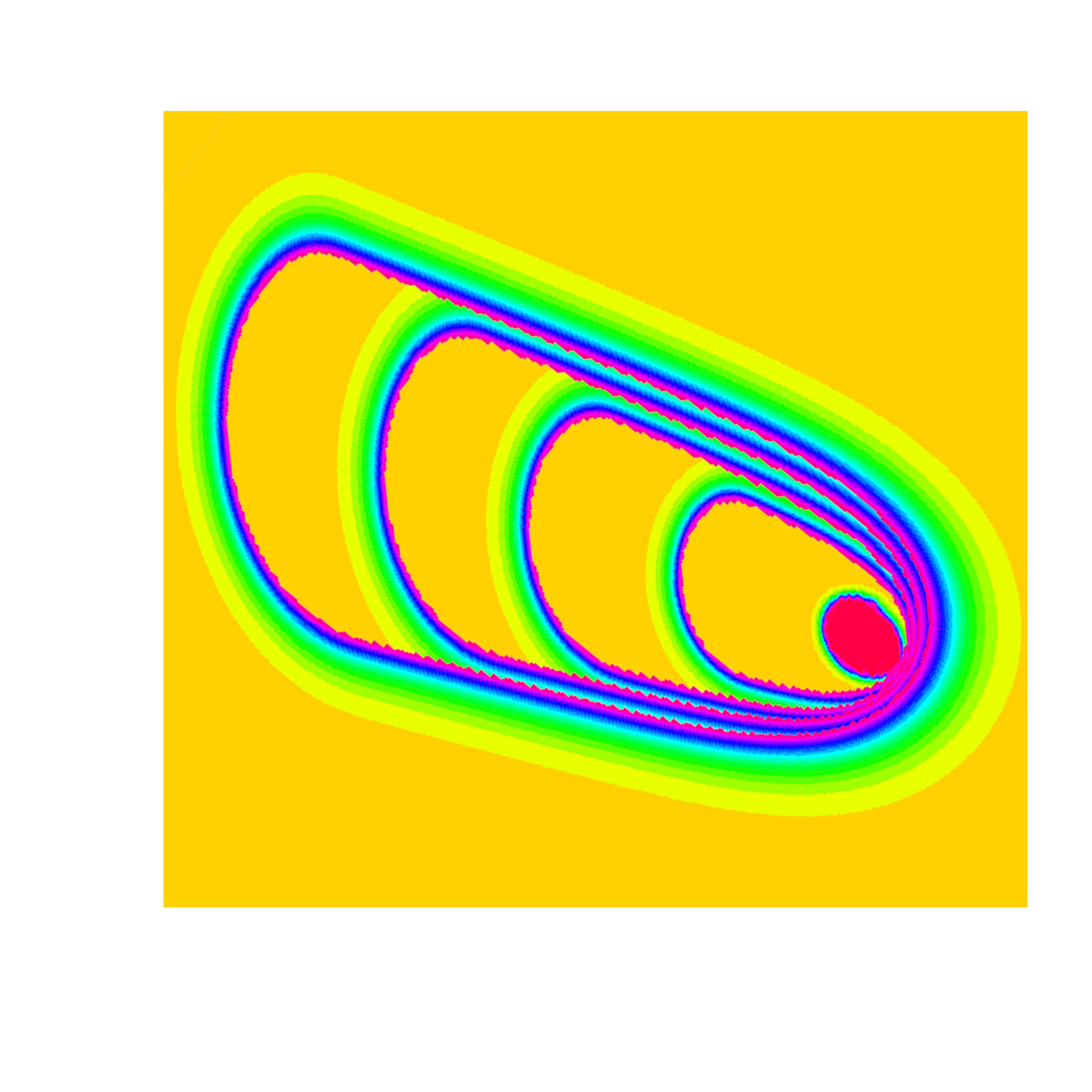}
    \end{minipage}
\caption{
Observed isochrones for spot fire E14 investigated in~\cite{CRUZ}, compared to
isochrones from our simulation taken at evenly spaced iterations of the scheme.
The red line is the fire front with temperature above ignition, and temperature
decreasing below ignition as the colours progress from blue to green, then
yellow.}
\label{FIG1}
\end{figure}

Numerical simulations were written in C++ with a standard finite-element
method, using the library FreeFEM++~ \cite{HECHT}.
Our simulations used the weak formulation of the PDE~\eqref{MAIN:EQ} for $u$ on a unit square
with Dirichlet zero boundary conditions.
For the convolution kernel $K(x,y)$ we used an approximation to the Dirac mass,
which allows close non-local interactions and facilitates fast computation.
For parameters, we used~$c = 10^{-3}$, $\Theta(y,t) = 1$, and $\beta(u) = 0$.
Figure \ref{FIG1} uses wind with~$\omega = (-1,0.4)$, otherwise~$\omega = 0$.
For the finite element method we used piecewise continuous elements and divided
the domain into a mesh with 10,000 individual elements.
Our simulations were run on a standard laptop without special processing power.

We remark that we did not make special effort to change parameters in~\eqref{MAIN:EQ}
drastically to fit pictures, but instead tried to make only natural, broad choices.
In a real-world scenario, we expect that parameters are at best approximated,
so fine-tuning them would be counterproductive in our view.

Specifically, in Figure~\ref{FIG1} we compare our numerical simulations with
one of the outdoor spot fire experiments conducted in Ridgewood by the CSIRO
in NSW Australia~\cite{CRUZ}. 
Named E14, the isochrones of the fire within 90 seconds from its ignition are
reproduced on the left, as recorded in~\cite{CRUZ}.
For fire E14, the authors in~\cite{CRUZ} note a change of wind direction by
approximately 40 degrees during its burning. Our simulation has a constant wind
direction.

On the right, we present the output of our numerical simulation using the
model in~\eqref{MAIN:EQ}.
The specifics of this numerical validation considered the fire emerging from a
point, with wind WNW and constant (24 degrees above the West direction).
There is a Dirichlet zero boundary condition imposed on the boundary of the
fire domain.
There is a slight difference in shape from our simulation compared to the
experiment, which is likely due to the volatile wind in the experiment.

\begin{figure}[h]
\centering

\includegraphics[width=.23\textwidth,height=0.2\textheight]{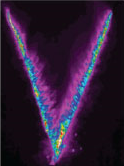}\hfil
\includegraphics[width=.23\textwidth,height=0.2\textheight]{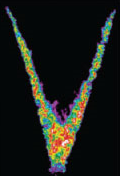}\hfil
\includegraphics[width=.23\textwidth,height=0.2\textheight]{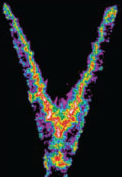}\hfil
\includegraphics[width=.23\textwidth,height=0.2\textheight]{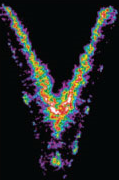}

\includegraphics[width=.23\textwidth,height=0.2\textheight,trim={60pt 0 60pt 0},clip]{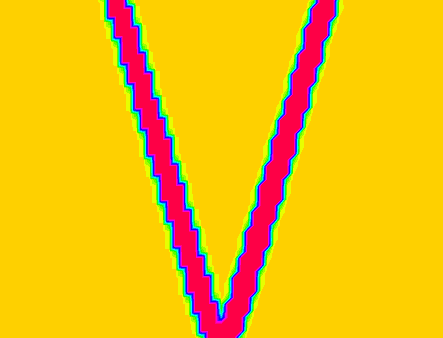}\hfil
\includegraphics[width=.23\textwidth,height=0.2\textheight,trim={60pt 0 60pt 0},clip]{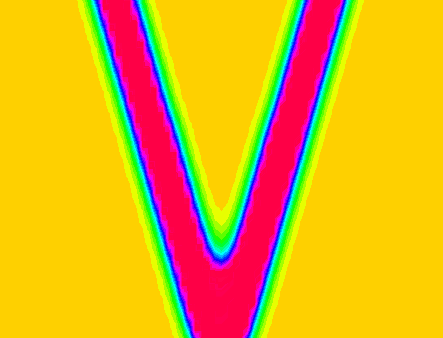}\hfil
\includegraphics[width=.23\textwidth,height=0.2\textheight,trim={60pt 0 60pt 0},clip]{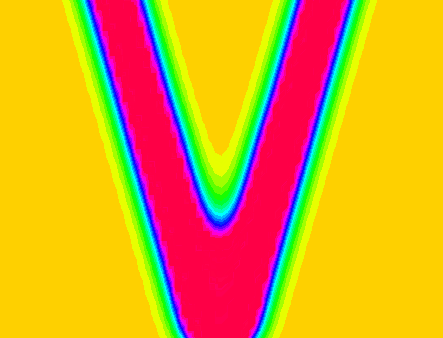}\hfil
\includegraphics[width=.23\textwidth,height=0.2\textheight,trim={60pt 0 60pt 0},clip]{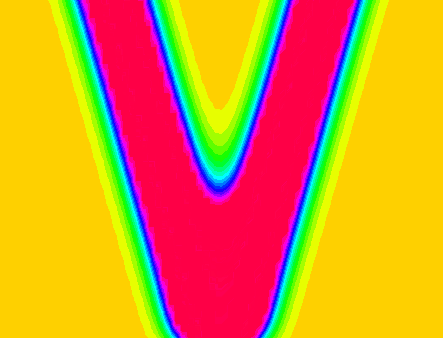}

\caption{Experiment CF-14 from~\cite{VIEGAS}, with a 30 degree cone angle. Frames taken from an infrared camera at thirty second intervals.
Simulation frames taken at evenly spaced iterations of the numerical simulation.}\label{FIG2}
\end{figure}

\begin{figure}[h]
\centering
\includegraphics[width=.31\textwidth]{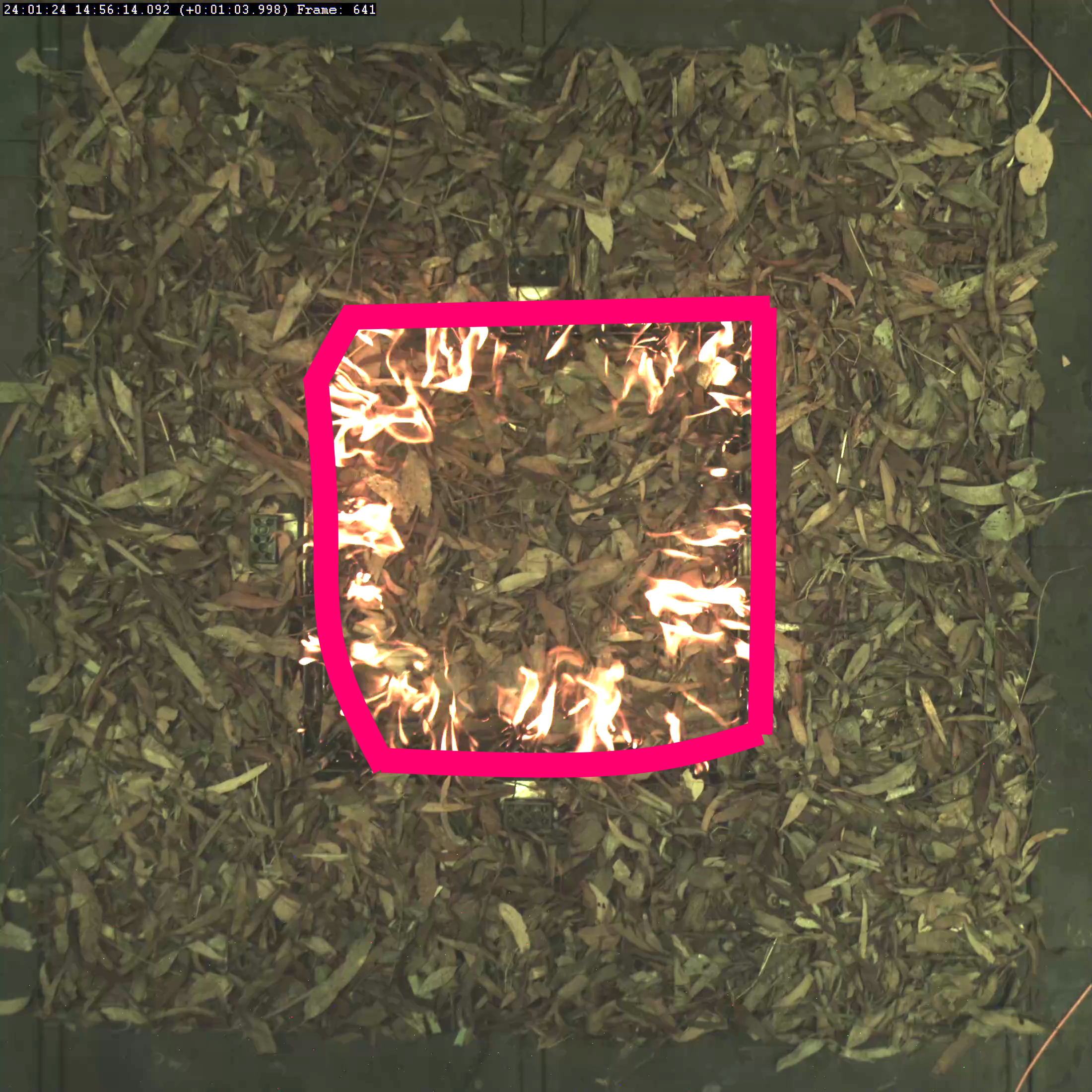}\hfil
\includegraphics[width=.31\textwidth]{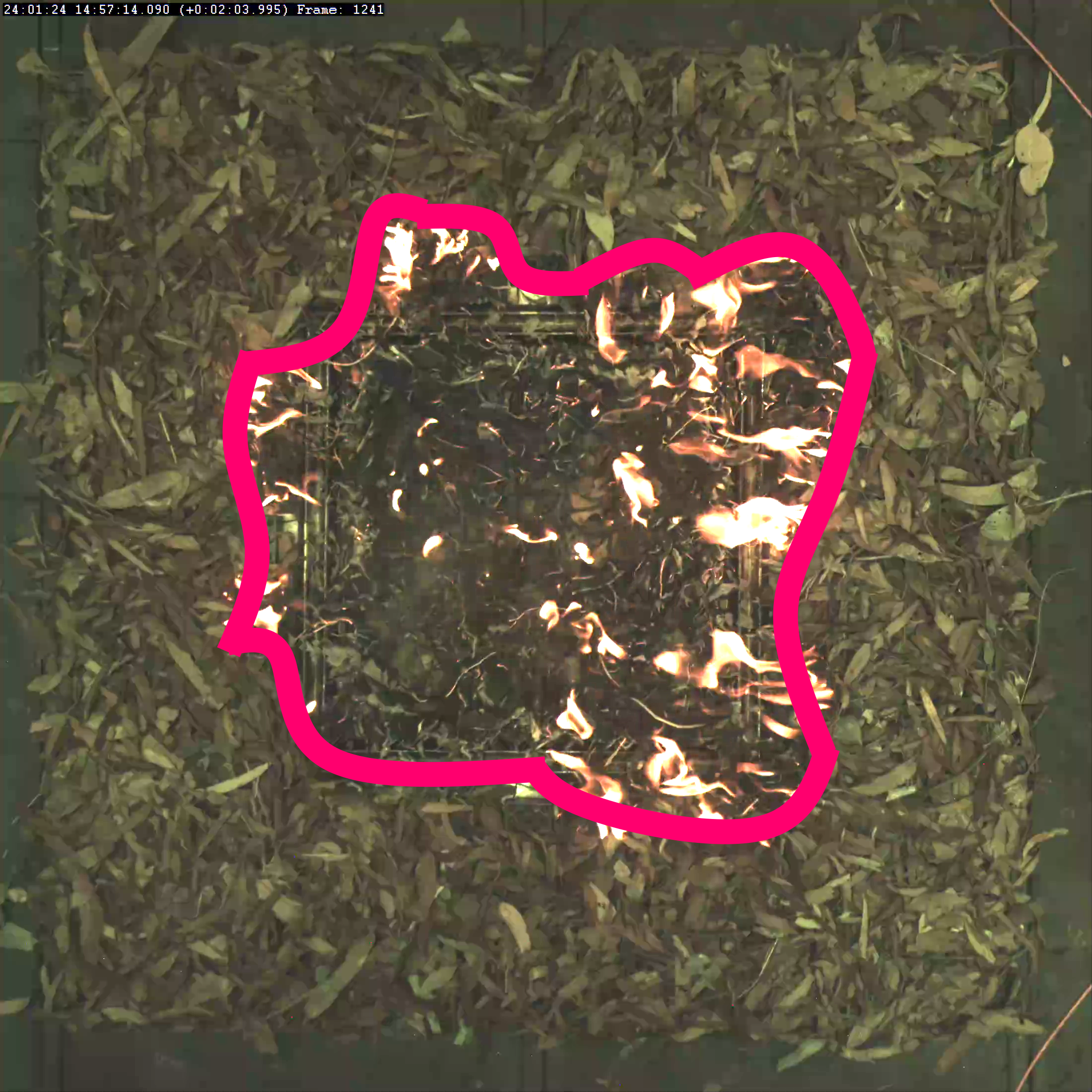}\hfil
\includegraphics[width=.31\textwidth]{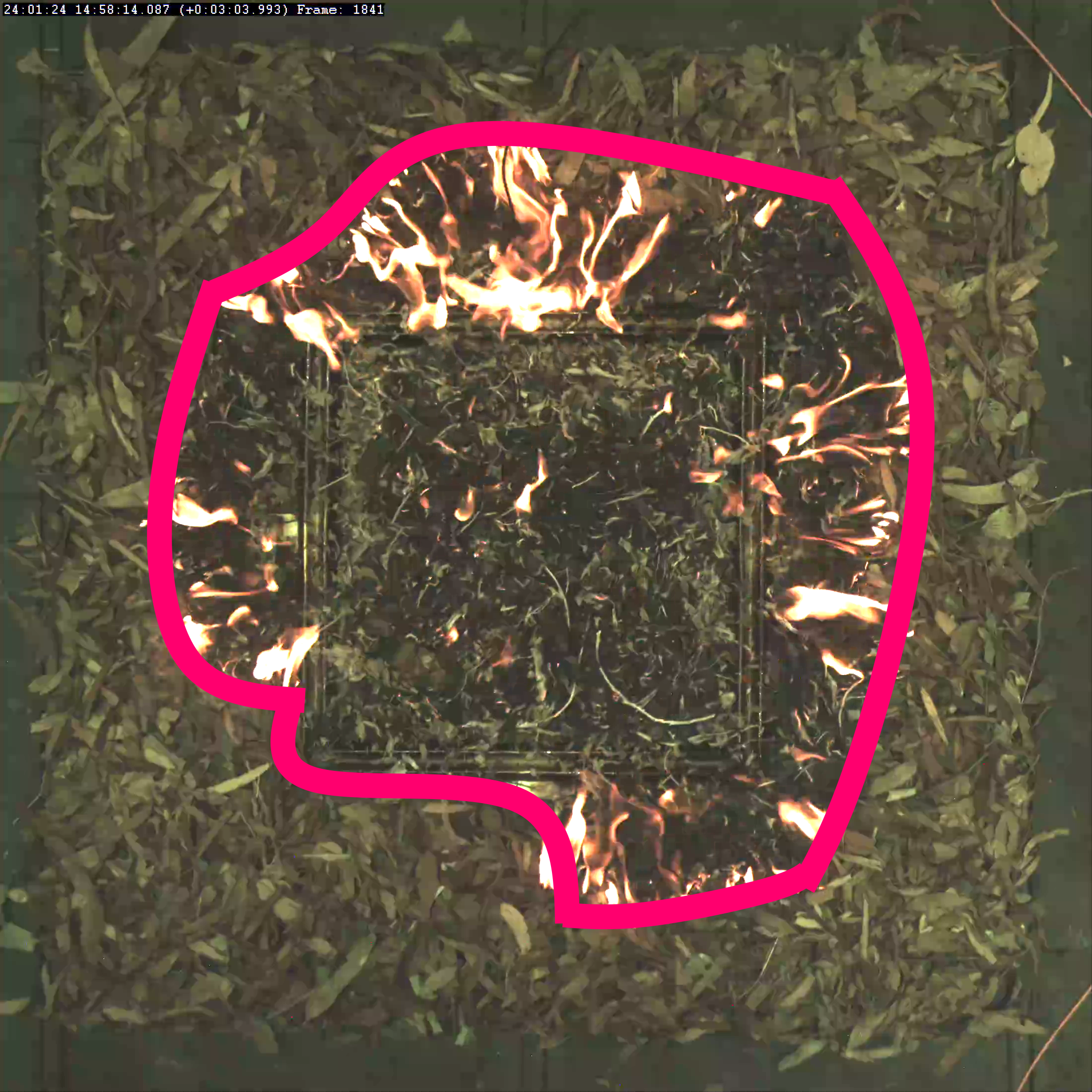}
\vspace{1mm}

\includegraphics[width=.31\textwidth]{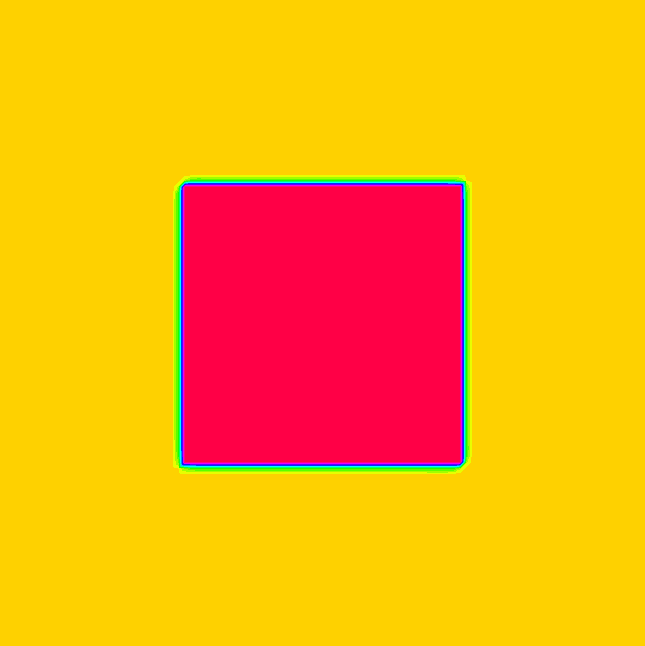}\hfil
\includegraphics[width=.31\textwidth]{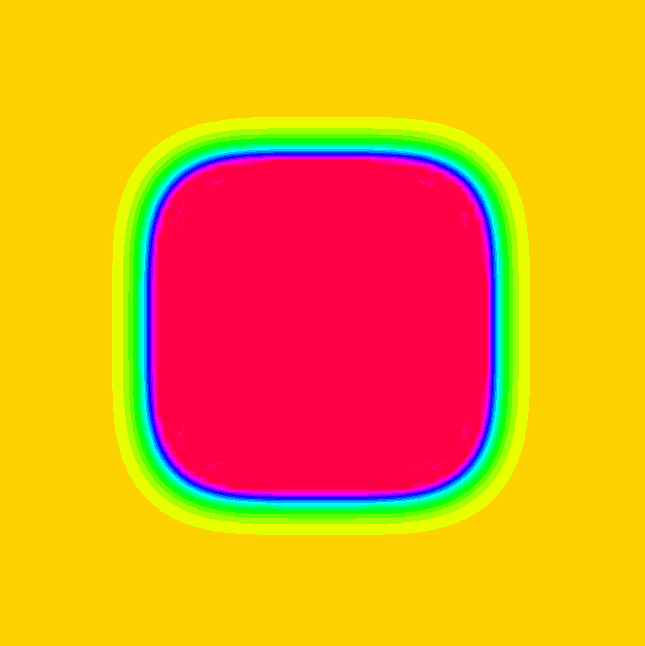}\hfil
\includegraphics[width=.31\textwidth]{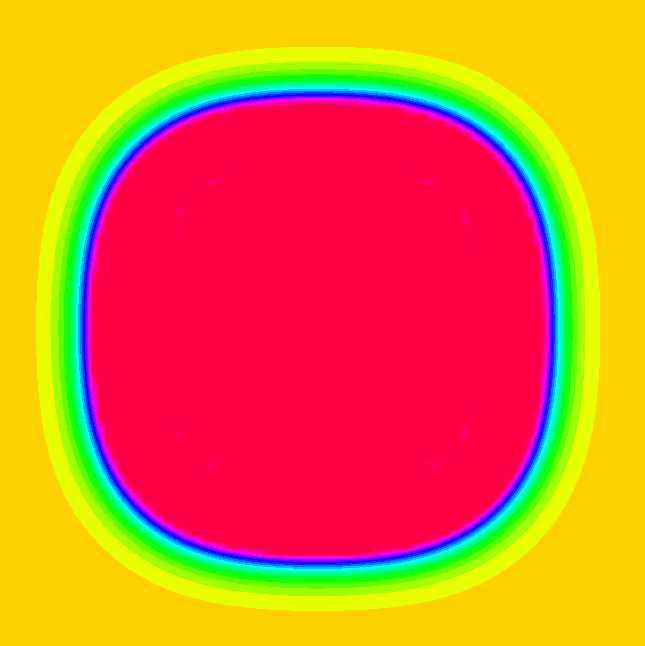}
\vspace{1mm}

\includegraphics[width=.31\textwidth]{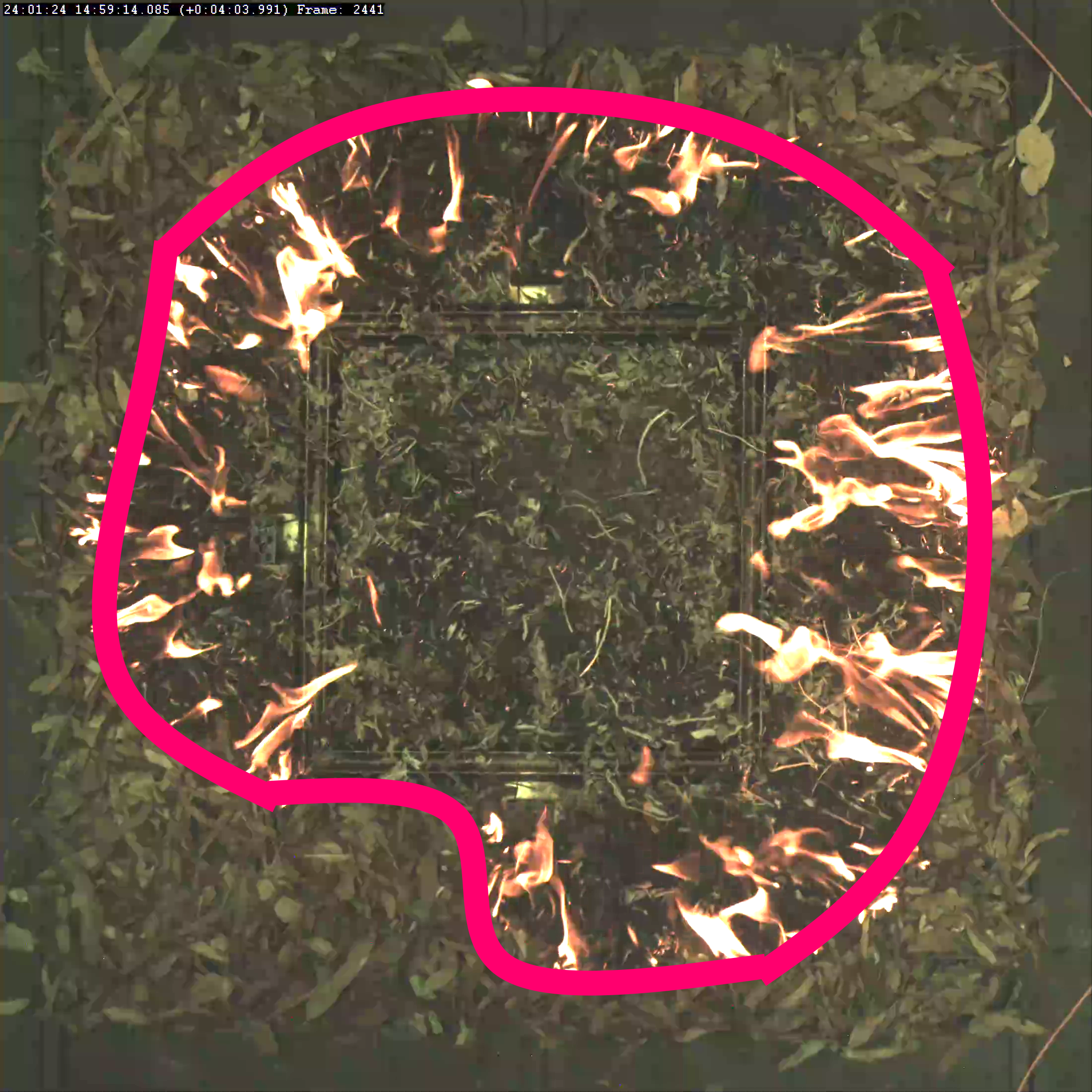}\hfil
\includegraphics[width=.31\textwidth]{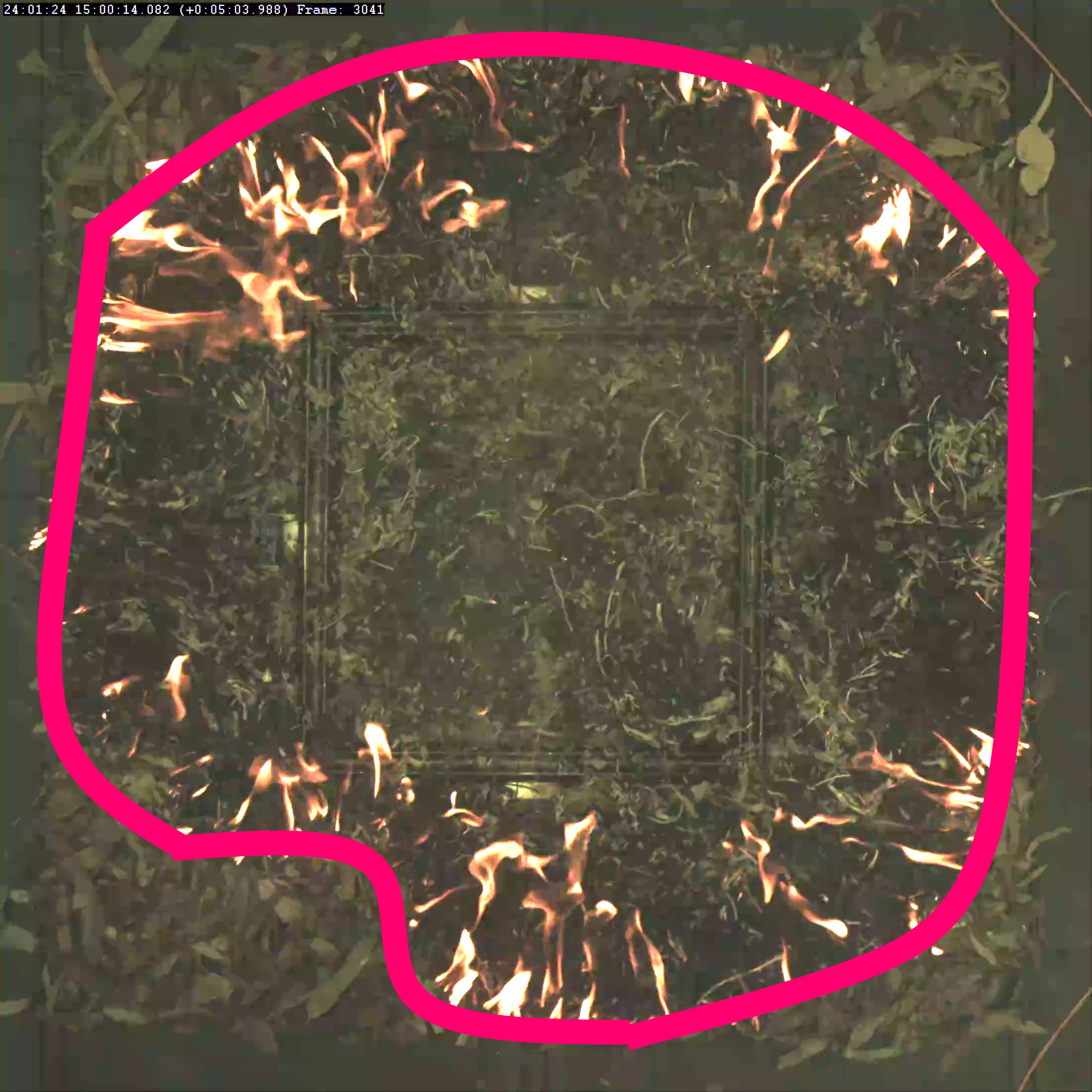}\hfil
\includegraphics[width=.31\textwidth]{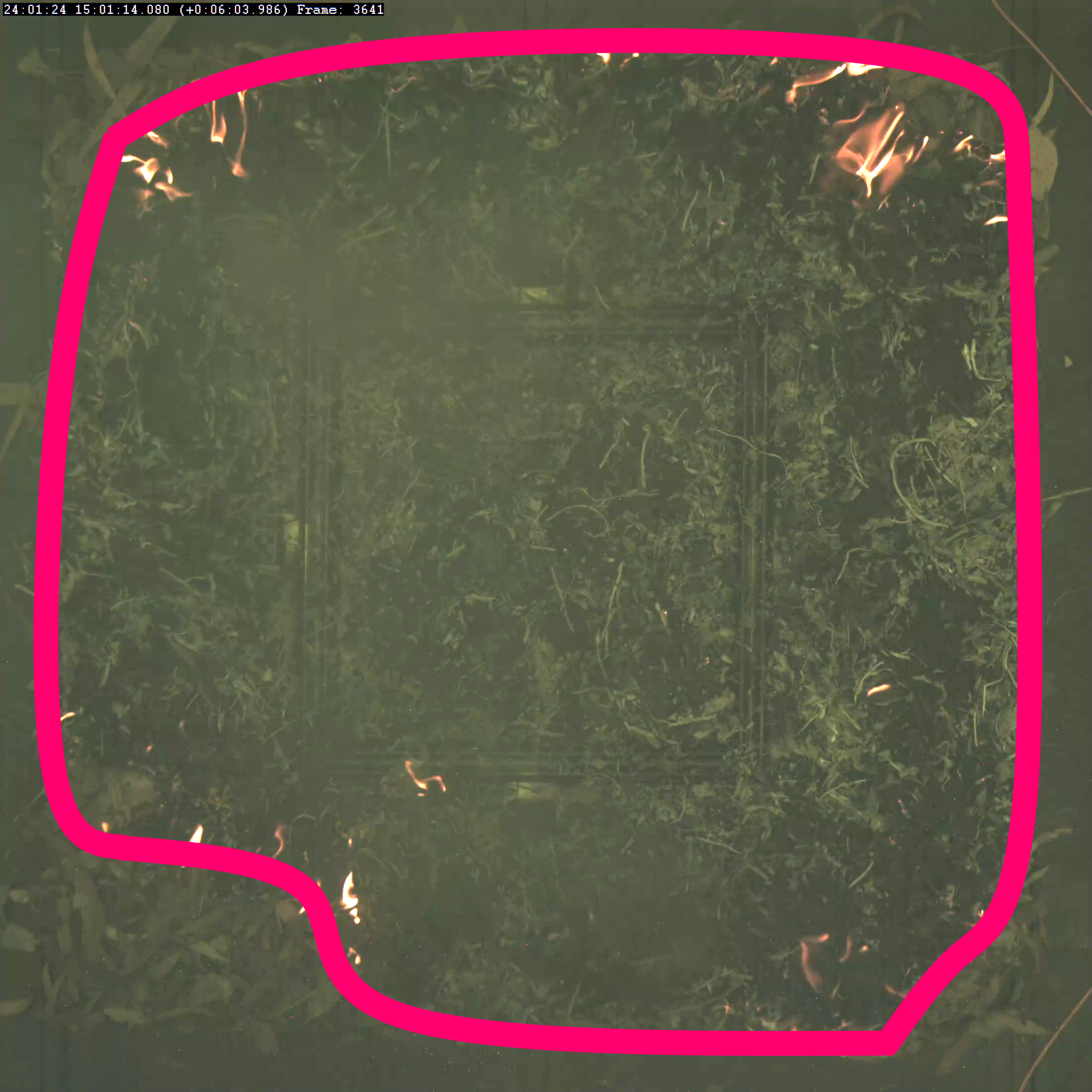}
\vspace{1mm}

\includegraphics[width=.31\textwidth]{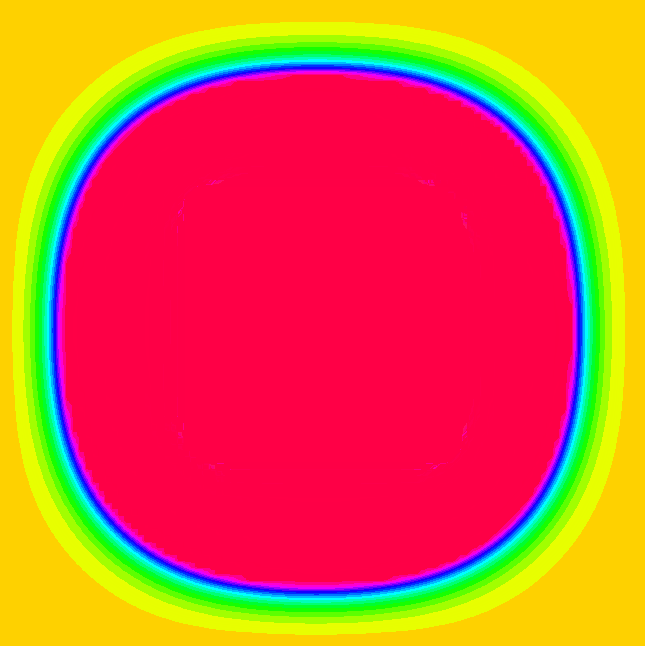}\hfil
\includegraphics[width=.31\textwidth]{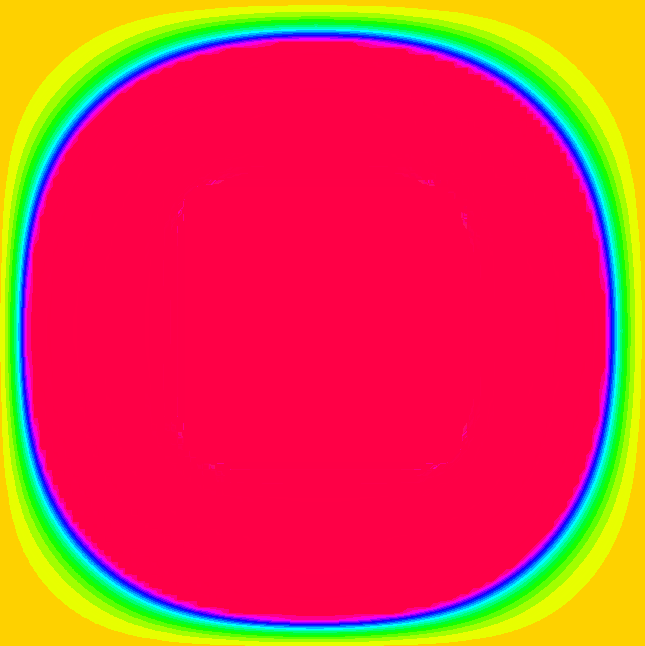}\hfil
\includegraphics[width=.31\textwidth]{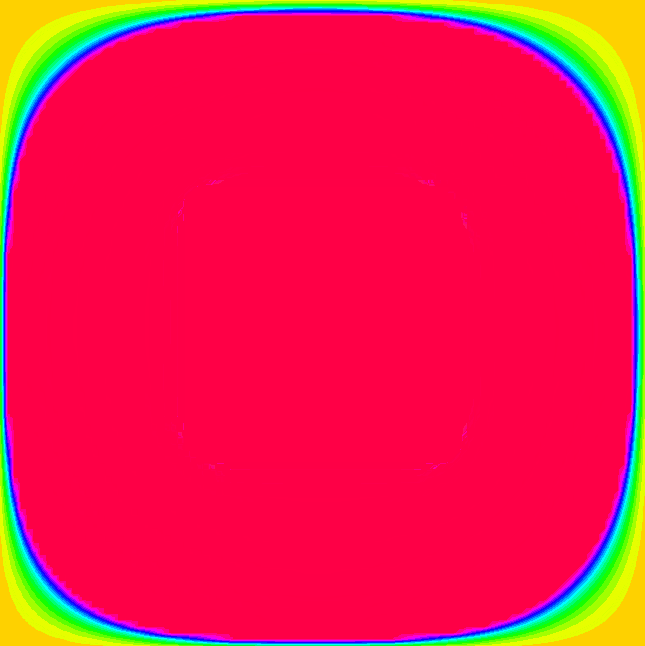}
\caption{Experimental vs. Simulation Frames. Experiment conducted at the Pyrotron CSIRO Canberra February 2024. Experimental frames taken at one-minute intervals. A magenta outline is traced around the approximate fire line. Simulation frames taken at evenly spaced iterations of the numerical simulation.}
\label{FigSquare}
\end{figure}

\begin{figure}
\centering
\includegraphics[width=.95\textwidth]{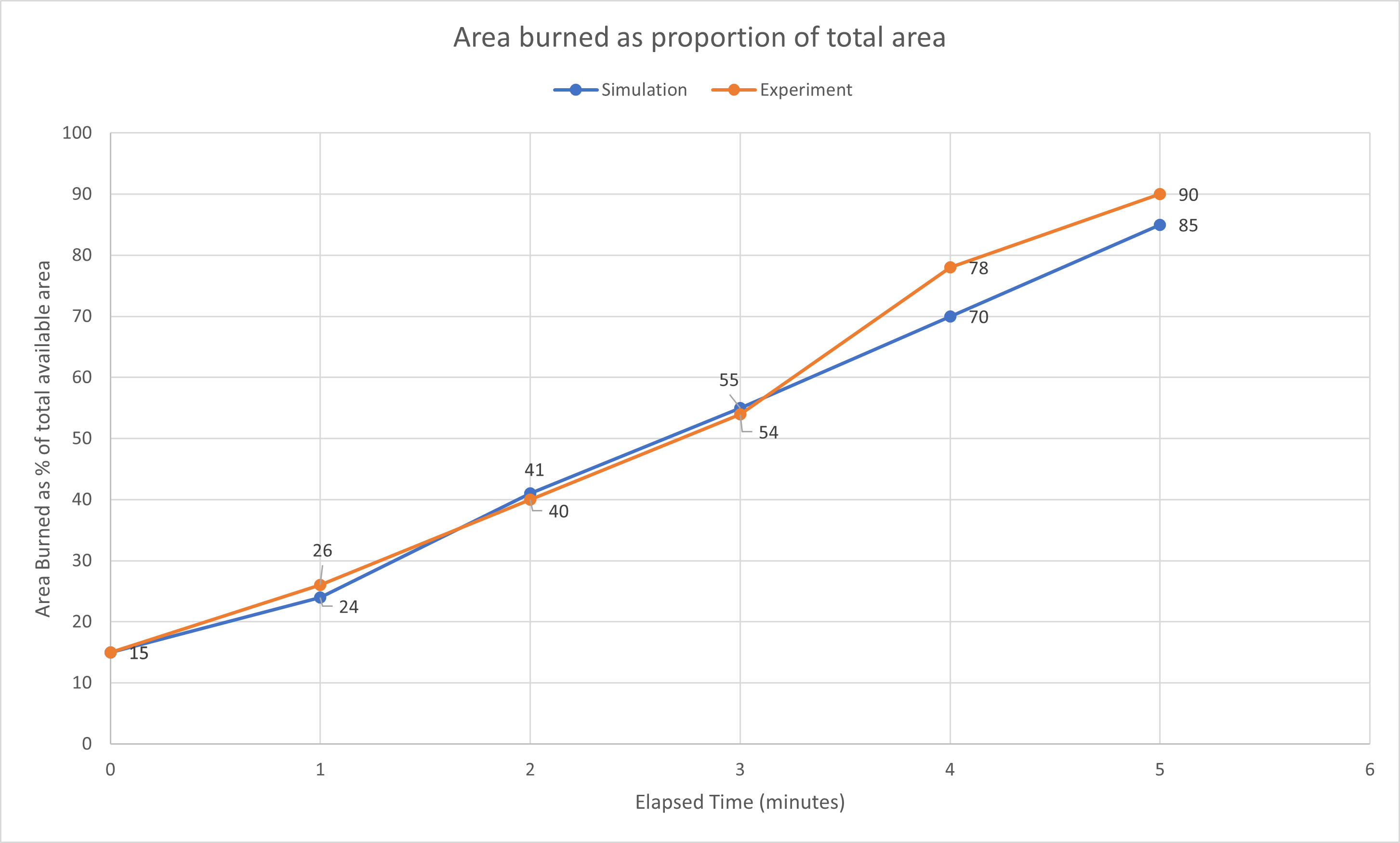}
\caption{Comparison of area within fire line for the experiment vs the simulation displayed in Figure \ref{FigSquare}.}\label{FigSquare2}
\end{figure}

In Figures~\ref{FIG2}, \ref{FigSquare}, and~\ref{FigSquare2} we compare our numerical simulations with experiments of fires run in labs.
Namely, Figure~\ref{FIG2} reproduces (on the top) the thermal images of the
experiment in~\cite{VIEGAS} of a V-shaped fire, to be confronted (on the
bottom) with our numerical simulations of the model in~\eqref{MAIN:EQ}.
This shows that our model works well even in the challenging case of fire-line
merging; the movement of the vertex of the fire is quite close to the observed
vertex, cf.~\cite{W15}.

In Figures~\ref{FigSquare} and~\ref{FigSquare2},
we report about an experiment recently taken place
at Pyrotron CSIRO Canberra, tracing the fire lines emanating from
an initial burning square. The photos taken in the lab are confronted
in Figures~\ref{FigSquare} with our simulations.

This comparison is retaken also in Figure~\ref{FigSquare2},
taking the approximated burned area as a measure of validation.
The matching is almost perfect after up to 3 minutes, and a small deviation
takes place at minute 4 and minute 5.  This small deviation has also to be
expected, since
our simulations are based on ideal Dirichlet conditions, while the experimental
fire may be affected by slightly approximated boundary conditions allowing for
a faster diffusion near the boundary.  Note also the coarseness of the fuel
used for the experimental fire.

\section{The hybrid geometric problem}

For constant ignition temperature~$\Theta_0$, one can consider the level set~$\{u=\Theta_0\}$ (or, more precisely, the boundary of~$\{u>\Theta_0\}$)
as the moving front of the fire. The velocity of this front in the external normal direction, in light of~\eqref{MAIN:EQ}, is 
\begin{equation}\label{KS-ev}
-cH+\frac{c\partial^2_\nu u}{|\nabla u|}+\frac1{|\nabla u|}\int_{\R^n}\big(u(y,t)-\Theta_0\big)_+\,K(x,y)\,dy+\left( \left(\omega
-
\frac{\beta(u)\nabla u}{|\nabla u|^\alpha}\right)\cdot\nabla u\right)_-\,,
\end{equation}
where~$H$ is the mean curvature of~$\partial\{u>\Theta_0\}$ and~$\nu$ the external unit normal.\medskip

The evolution equation in~\eqref{KS-ev} has an interesting structure
since it mixes a geometric term of curvature type with a forcing term arising from the global pattern of the solution of~\eqref{MAIN:EQ}. This ``hybrid'' structure reflects the physical property that the moving front of the bushfire is not modeled solely by its geometry, since it is influenced by all the structural and environmental parameters, such as temperature, wind, and convection.

Similar, but structurally different, equations have been obtained, or imposed a-priori, in other bushfire models, see e.g.~\cite[equation~(2.6)]{MR4439510}, in which the outward normal velocity of a fire line in the plane was set to be equal (up to normalizing constants) to a constant (inherited by free boundary problems of Stefan type and accounting for oxygen gradient), minus 
the curvature (which was proposed in~\cite{MARK}
as a corrective term to match the experimental data
and acts as a smoothing term on perturbations of the fire line, see e.g.~\cite[Figure~5]{MR4439510}), plus a dynamical term requiring the solution of an advection-diffusion for the oxygen concentration.
\medskip

Interestingly, in our model, the normal velocity of the front is not imposed a-priori, but
follows from the constitutive equation~\eqref{MAIN:EQ}.
To see how formally deduce~\eqref{KS-ev} from~\eqref{MAIN:EQ} when~$\Theta=\Theta_0$, let us
consider the level set~$\{u=\Theta_0\}$ and suppose, up to a rigid motion, that it passes through the origin of~$\R^n$ at time~$t=0$, and that, in a neighborhood of the origin, the burning region~$\{u>\Theta_0\}$ corresponds to~$\{x_n>\gamma(x',t)\}$ for some~$\gamma:\R^{n-1}\times\R\to\R$, with~$\gamma(0,0)=0$ and
\begin{equation}\label{GRA}
\nabla_{x'} \gamma(0,0)=0,\end{equation}
where we adopt the notation~$ x' = (x_1,\ldots,x_{n-1}) $.

This rigid motion is convenient from the computational point of view, since it allows us to compute the mean curvature~$H$ at the origin\footnote{In our notation, the mean curvature is the tangential divergence of the normal vector field; as such it is not normalized by a factor~$\frac1{n-1}$.} as
$$\left.{\rm div}_{x'}\left(\frac{\nabla_{x'} \gamma}{\sqrt{1+|\nabla_{x'} \gamma|^2}}\right)\right|_{(0,0)}=\Delta_{x'}\gamma(0,0),$$ see e.g.~\cite[Theorem~2.3.7]{ZZL}.

Then, differentiating the relation~$u(x',\gamma(x',t),t)=\Theta_0$ with respect to the space variables, we obtain that, for each~$j$, $k\in\{1,\dots,n-1\}$,
\begin{eqnarray*}
&& \partial_{x_j} u(x',\gamma(x',t),t)+\partial_{x_n} u(x',\gamma(x',t), t)\,\partial_{x_j}\gamma(x',t)=0\\
{\mbox{and }}&&
\partial_{x_j x_k}^2 u(x',\gamma(x',t),t)+
\partial_{x_j x_n}^2 u(x',\gamma(x',t),t)\,\partial_{x_k}\gamma(x',t)+\partial_{x_k x_n}^2 u(x',\gamma(x',t),t)\,\partial_{x_j}\gamma(x',t)\\&&\qquad
+\partial_{x_n}^2 u(x',\gamma(x',t), t)\,\partial_{x_j}\gamma(x',t)\,\partial_{x_k}\gamma(x',t)+\partial_{x_n} u(x',\gamma(x',t), t)\,\partial^2_{x_j x_k}\gamma(x',t)
=0.
\end{eqnarray*}
We now evaluate the last equation at~$x'=0$ and~$t=0$, using~\eqref{GRA}, and we conclude that, for each~$j$, $k\in\{1,\dots,n-1\}$,
$$ 
\partial_{x_j x_k}^2 u(0,0)+\partial_{x_n} u(0,0)\,\partial^2_{x_j x_k}\gamma(0,0)
=0$$
and therefore, choosing~$j=k$ and summing up,
\begin{equation}\label{12w}\begin{split}
0&=\sum_{k=1}^{n-1}\partial_{x_k}^2 u(0,0)+\partial_{x_n} u(0,0)\,\Delta_{x'}\gamma(0,0)\\
&=\Delta u(0,0)-\partial_{x_n}^2 u(0,0)+|\nabla u(0,0)|\,H(0,0).
\end{split}\end{equation}

Besides, differentiating the relation~$u(x',\gamma(x',t),t)=\Theta_0$ with respect to the time variables, we obtain that
\begin{eqnarray*}
\partial_{x_n}u(x',\gamma(x',t),t)\,\partial_t\gamma(x',t)
+\partial_t u(x',\gamma(x',t),t)=0
\end{eqnarray*}
and therefore
$$ |\nabla u(0,0)|\,\partial_t\gamma(0,0)+\partial_t u(0,0)=0.$$

{F}rom this, \eqref{MAIN:EQ} and~\eqref{12w}, omitting the evaluation at~$(0,0)$ for short, we arrive at
\begin{eqnarray*}
-\partial_t\gamma&=&\frac{\partial_t u}{|\nabla u|}\\
&=&\frac{c\Delta u}{|\nabla u|}+\frac1{|\nabla u|}\int_{\R^n}\big(u(y,t)-\Theta(y,t)\big)_+\,K(x,y)\,dy+\left( \left(\omega-\frac{\beta(u)\nabla u}{|\nabla u|^\alpha}\right)\cdot\frac{\nabla u}{|\nabla u|}\right)_-\\
&=&-cH+\frac{c\partial_{x_n}^2 u}{|\nabla u|}
+\frac1{|\nabla u|}\int_{\R^n}\big(u(y,t)-\Theta(y,t)\big)_+\,K(x,y)\,dy\\&&\qquad\qquad+\left( \left(\omega-\frac{\beta(u)\nabla u}{|\nabla u|^\alpha}\right)\cdot\frac{\nabla u}{|\nabla u|}\right)_-\,,
\end{eqnarray*}
which is~\eqref{KS-ev}.

\section{Conclusions}

We have put forth a model for bushfires relying on partial differential equations. In its simplest form, the model takes the form
\begin{equation*}\begin{split}
\partial_t u&=c\Delta u+\int_{\R^n}\big(u(y,t)-\Theta(y,t)\big)_+\,K(x,y)\,dy+\left( \left(\omega-\frac{\beta(u)\nabla u}{|\nabla u|^\alpha}\right)\cdot\nabla u\right)_-\,,\end{split}
\end{equation*}
where~$u$ represents the temperature of the environment,
$c$ the environmental diffusion coefficient, $\Theta$ the effective ignition temperature, $\omega$ a parameter accounting for wind, and~$\beta$ a function modulating a convective effect.

The front of the fire can be modeled by the level set of the temperature function corresponding to the ignition temperature.
The velocity of the front in the normal direction
when the ignition temperature~$ \Theta_0 $ is constant takes the form
$$ -cH+\frac{c\partial^2_\nu u}{|\nabla u|}
+\frac1{|\nabla u|}\int_{\R^n}\big(u(y,t)-\Theta_0\big)_+\,K(x,y)\,dy+\left( \left(\omega-\frac{\beta(u)\nabla u}{|\nabla u|^\alpha}\right)\cdot\nabla u\right)_-\,,$$
where~$H$ is the mean curvature of the burning region.

Specific features of this setting include:
\begin{itemize}
\item The basic model can be easily generalized to comprise additional phenomena,
including a memory effect to identify regions which have been already burned, variations of the environment in space and time, different types of heat diffusions along the ground and in the air, etc.
\item Albeit non-standard, the partial differential equation is endowed with a solid mathematical theory,
\item The description of the geometric front combines a smoothing effect, encoded by the mean curvature term (with a positive sign,
allowing for the descriptions of moving fronts starting, for instance,
by domains with corners),
\item The model can be efficiently implemented through numerical simulations,
\item The simulations run in real-time and are in agreement both with the data collected in supervised lab experiments and with those coming from real bushfire events. 
\end{itemize}

\vfill

\end{document}